\documentclass{article}
\usepackage{spconf,amsmath,graphicx}
\usepackage{amssymb}
\usepackage{amsfonts}
\usepackage{physics}

\newcommand*{\QEDB}{\hfill\ensuremath{\blacksquare}}%

\newcommand{\hide}[1]{}

\usepackage{algorithm} 
\usepackage{algpseudocode} 

\usepackage{subcaption}
\usepackage{caption}
\usepackage{xcolor}
\def\mz{\mu_z}
\def\mt{\mu_\theta}
\title{A COMMUNICATION EFFICIENT QUASI-NEWTON METHOD FOR LARGE-SCALE DISTRIBUTED MULTI-AGENT OPTIMIZATION}
\name{Yichuan Li$^1$,  Petros G. Voulgaris$^{2}$, and Nikolaos M. Freris$^3$\thanks{This work was supported by the Ministry of Science and Technology
of China under grant 2019YFB2102200 and the Anhui Dept. of Science and
Technology under grant 201903a05020049. Corresponding author is N. Freris.}}
\address{$^1$Coordinated Science Laboratory, University of Illinois at Urbana-Champaign\\$^{2}$Department of Mechanical Engineering, University of Nevada, Reno\\$^{3}$School of Computer Science, University of Science and Technology of China}

\begin{document}
%
\maketitle
\begin{abstract}
We propose a communication efficient quasi-Newton method for large-scale multi-agent convex composite optimization. We assume the setting of a network of agents that cooperatively solve a global minimization problem with strongly convex local cost functions augmented with a non-smooth convex regularizer. By introducing consensus variables, we obtain a block-diagonal Hessian and thus eliminate the need for additional communication when approximating the objective curvature information. Moreover, we reduce computational costs of existing primal-dual quasi-Newton methods from $\mathcal{O}(d^3)$ to $\mathcal{O}(cd)$ by storing $c$ pairs of vectors of dimension $d$. An asynchronous implementation is presented that removes the need for coordination. Global linear convergence rate in expectation is established, and we demonstrate the merit of our algorithm numerically with real datasets. 
\end{abstract}


\section{Introduction}
\label{sec1}
Distributed multi-agent optimization \cite{distOpt} has found numerous applications in a range of fields, such as compressed sensing \cite{compressedsensing}, distributed learning \cite{learning}, and sensor networks \cite{sensor}. A canonical problem can be formulated as follows: 
\begin{gather}
    \underset{\hat{w}\in\mathbb{R}^d}{\text{minimize}} \,\bigg\{\sum_{i=1}^m f_i(\hat{w})+g(\hat{w}) \bigg\}, \label{prob1}
\end{gather} 
where $f_i(\cdot)$ is a strongly convex local cost function pertaining to agent $i$
and $g(\cdot)$ is a convex but possibly nonsmooth local regularizer (for example, the $\ell_1$-norm that aims to promote sparsity in the solution). We assume the setting of a network of agents who cooperatively solve (\ref{prob1}) for the common decision variable $\hat{w}$ by exchanging messages with their neighbors. First order methods \cite{subgrad}\nocite{pgextra}-\cite{p2d2} have been popular choices for tackling problem (\ref{prob1}) due to their simplicity and economical computational costs. However, such methods typically feature slow convergence rate,\hide{ i.e., may require a large number of iterations to reach a prescribed accuracy.} whence a remedy is to resort to higher-order information to accelerate convergence. Newton methods use the Hessian of the objective to scale the gradient so as to accelerate the convergence speed. However, Newton methods suffer from two main drawbacks that hinder their broader applicability in distributed optimization: (i) they are limited to low dimensional problems since, at each iteration, a linear system of equations has to be solved which incurs $\mathcal{O}(d^3)$ computational costs, (ii) in multi-agent networks, the Newton step for each agent depends on its neighbors. In other words, computing the Newton step would involve multiple rounds of communication as in \cite{nn}\nocite{dsop}\nocite{asy-nn}-\cite{dnadmm}. \\
\textbf{Contributions}: (i) We propose a communication efficient quasi-Newton method based on the Alternating Direction Method of Multipliers (ADMM) \cite{boyd}. We decouple agents from their neighbors through the use of consensus variables. This achieves a block-diagonal Hessian so that quasi-Newton steps can be computed locally \emph{without additional communication among agents}. (ii) By storing $c$ number of vector pairs, we reduce the computational costs of existing primal-dual quasi-Newton methods from $\mathcal{O}(d^3)$ to $\mathcal{O}(cd)$. (iii) We present an asynchronous implementation scheme that only involves a subset of agents at each iteration, thus removing the need for network-wide coordination. (iv) We establish global linear convergence rate in expectation without backtracking, and demonstrate its merits with numerical experiments on real datasets.

\section{Preliminaries}
\label{sec2}
\subsection{Problem formulation}
The network of agents is captured by an undirected graph $\mathcal{G}= (\mathcal{V},\mathcal{E})$, where $\mathcal{V}=\{1,\dots,m\}$ is the vertex set and $\mathcal{E}\subseteq \mathcal{V}\times \mathcal{V}$ is the edge set, i.e., $(i,j)\in \mathcal{E}$ if and only if agent $i$ can communicate with agent $j$. The total number of communication links is denoted by $n=\abs{\mathcal{E}}$ and the neighborhood of agent $i$ is represented as $\mathcal{N}_i:=\{j\in \mathcal{V}:(i,j)\in \mathcal{E}\}$. We further define the source and destination matrices $\hat{A}_s,\hat{A}_d\in \mathbb{R}^{n\times m}$ as follows: we let the $k$-th edge be represented by the $k$-th row of $\hat{A}_s$ and $\hat{A}_d$, i.e., $[\hat{A}_s]_{ki}=[\hat{A}_d]_{kj}=1$, while all other entries in the $k$-th row are zero. Using the above definition, we introduce local copies of the decision variable, i.e, $w_i$ held by agent $i$, and reformulate problem (\ref{prob1}) to the consensus setting as follows:
\begin{gather}
        \underset{w_i,\,\theta,\,z_{ij}\,\in \mathbb{R}^d}{\mathrm{minimize}}\left\{\sum_{i=1}^m f_i(w_i)+g(\theta)\right\}\nonumber\\
    \text{s.t.} \quad 
    w_i=z_{ij}=w_j\quad \forall\, i\,\,\text{and}\, \, j\in \mathcal{N}_i,\nonumber\\
    w_l = \theta \quad \text{for some }\,l\in[m].\label{prob2}
\end{gather}
We separate the decision variable corresponding to the smooth part and the nonsmooth part of the objective by enforcing an equality constraint: $w_l=\theta$, for arbitrarily chosen $l\in[m]$. Moreover, the consensus variable $\{z_{ij}\}$ decouples $w_i$ from $w_j$ which is \textit{crucial} to achieve a communication-efficient implementation. Through the use of $\{z_{ij}\}$, we show that the curvature of agent $i$ does not depend on its neighbors, whence only local information is needed to construct estimates. See Section 3 for detailed discussion. When the network graph is connected, problem (\ref{prob2}) is equivalent to problem (\ref{prob1}) in the sense that the optimal solution coincides, i.e., $\hat{w}^\star=w_i^\star=\theta^\star=z_{ij}^\star$ for all $i\in [m]$ and $j\in \mathcal{N}_i$. By stacking $w_i,z_{ij}$ into column vectors $w\in \mathbb{R}^{md}, z\in \mathbb{R}^{nd}$ respectively, we define $F(w)=\sum_{i=1}^m f_i(w_i)$, whence problem (\ref{prob2}) can be compactly expressed as:
\begin{gather}
    \underset{\theta\in\mathbb{R}^d,w\in \mathbb{R}^{md},z\in\mathbb{R}^{nd}}{\mathrm{minimize}} \bigg\{F(w)+g(\theta)\bigg\} \nonumber\\
    \text{s.t.}\quad Aw=\begin{bmatrix}\hat{A}_s \otimes I_d \nonumber\\
    \hat{A}_d\otimes I_d\end{bmatrix}w = \begin{bmatrix}I_{nd}\nonumber\\I_{nd} \end{bmatrix}z=Bz, \nonumber\\
    S^\top w =\theta, \label{prob3}
\end{gather}
where $\otimes$ denotes the Kronecker product, $A\in \mathbb{R}^{2nd\times md}$ is formed by using the aforementioned source and destination matrices, and $S :=(s_l\otimes I_d)$ is formed by using the vector $s_l\in\mathbb{R}^m$ whose entries are zero except the $l$-th entry being 1. 
\subsection{Introduction to L-BFGS}
Quasi-Newton methods \cite{numerical} constitute a popular substitute for Newton methods when the dimension $d$ makes solving for the Newton updates computationally burdensome. In contrast to the Newton method, quasi-Newton methods do not directly invoke the Hessian but instead seek to estimate the curvature through finite differences of iterate and gradient values. In subsequent presentation, we focus on the quasi-Newton method proposed by Broyden, Fletcher, Goldfarb, and Shannon (BFGS) \cite{bfgs1}-\cite{bfgs2}. Specifically, we define the iterate difference $s^t:=w^{t+1}-w^t$ and the gradient difference $q^t:=\nabla f^{t+1}-\nabla f^t$. A Hessian inverse approximation scheme is iteratively obtained as follows:
\begin{gather}
    \left(\widehat{H}^{t+1}\right)^{-1} = \left(V^t\right)^\top \left(\widehat{H}^t\right)^{-1} V^t+\rho^t s^t \left(s^t\right)^\top, \label{bfgs}
\end{gather}
where $\rho^t = 1/\langle q^t,s^t\rangle$, and $V^t = I-\rho^t q^t (s^t)^\top$. In the BFGS scheme \cite{bfgsadmm}, the approximation $(\widehat{H}^t)^{-1}$ is stored and thus the update can be directly computed as $r^t = -(\widehat{H}^t)^{-1} \nabla f^t$ \emph{without solving} a linear system as in the Newton method: this reduces the computation costs from $\mathcal{O}(d^3)$ to $\mathcal{O}(d^2)$ for general problems. Nevertheless, (\ref{bfgs}) approximates $(\widehat{H}^t)^{-1}$ as a function of the initialization $(\widehat{H}^{0})^{-1}$ and $\{s^k,q^k\}_{k=0}^{t-1}$, i.e., the entire history of the iterates. Since early iterates tend to carry less information on the current curvature, Limited-memory BFGS (L-BFGS) \cite{lbfgs} was proposed to store only $\{s^k,q^k\}_{k=t-c}^{t-1}$ to estimate the update direction using Algorithm 1 \cite{numerical}. Note that L-BFGS operates solely on vectors, thus it not only reduces the storage costs from $\mathcal{O}(d^2)$ to $\mathcal{O}(cd)$, but also only requires $\mathcal{O}(cd)$ computation.

\begin{algorithm}[t]
	\caption{Two-Loop recursion of L-BFGS} 
	\textbf{Input}: $h=\nabla f^t$, $\{s^k,q^k\}_{k=t-c}^{t-1}$
	\begin{algorithmic}[1]
		\For {$k=t-1,\dots,t-c$}
		\State $\alpha^k\leftarrow\rho^k \langle s^k,h \rangle$
        \State $h \leftarrow h-\alpha^k q^k$
		\EndFor 
		\State $r \leftarrow (\widehat{H}^{t,0})^{-1}h$
		\For {$k=t-c,\dots,t-1$}
		\State $\beta\leftarrow\rho^k\langle q^k, r\rangle$ 
		\State $r \leftarrow r + (\alpha^k-\beta)s^k$
		\EndFor
	\end{algorithmic} 
	\textbf{Output}: $r = (\widehat{H}^{t,c})^{-1}\nabla f^t$
\end{algorithm}

\section{Proposed Method}
\label{sec3} 

ADMM solves (\ref{prob3}) and equivalently (\ref{prob1}) by operating on the Augmented Lagrangian defined as: 
\begin{gather*}
    \mathcal{L}(w,\theta,z;y,\lambda) = F(w)+g(\theta)+y^\top (Aw-Bz)\\+\lambda^\top (S^\top w-\theta)+\tfrac{\mz}{2}\norm{Aw-Bz}^2+\tfrac{\mt}{2}\norm{S^\top w-\theta}^2.
\end{gather*}
At each iteration, ADMM sequentially minimizes $\mathcal{L}(\cdot)$ with respect to primal variables $(w,\theta,z)$ and then performs gradient ascent on dual variables $(y,\lambda)$. Specifically, 
\begin{subequations}
\begin{align}
    w^{t+1} &= \underset{w}{\text{argmin}}\,\, \mathcal{L}(w,\theta^t,z^t;y^t,\lambda^t) ,\label{admm1}\\
    \theta^{t+1} & = \underset{\theta}{\text{argmin}}\,\,\mathcal{L}(w^{t+1},\theta,z^t;y^t,\lambda^t),\label{admm2}\\
    z^{t+1} & = \underset{z}{\text{argmin}}\,\,\mathcal{L}(w^{t+1},\theta^{t+1},z;y^t,\lambda^t),\label{admm3}\\
    y^{t+1} & = y^t + \mz (Aw^{t+1}-Bz^{t+1}) ,\label{admm4}\\
    \lambda^{t+1} & = \lambda^t +\mt(S^\top w^{t+1}-\theta^{t+1}).\label{admm5}
\end{align}
\end{subequations}
Note that above updates fall into the category of 3-block ADMM which is not guaranteed to converge for arbitrary $\mz,\mt>0$ \cite{3block}. Moreover, step (\ref{admm1}) requires the solution of a sub-optimization problem: this typically involves multiple inner-loops and thus induces heavy computational burden on agents. We propose to approximate step (\ref{admm1}) by performing the following \emph{one-step} update:
\begin{gather}
    w^{t+1} = w^t - (\widehat{H}^t)^{-1}\nabla_w \mathcal{L}(w^t,\theta^t,z^t;y^t,\lambda^t), \label{vaniprimal}
\end{gather}
where $(\widehat{H}^t)^{-1}$ is obtained by using L-BFGS (Algorithm 1) for an appropriate choice of $\{s^k,q^k\}$ pairs. Step (\ref{admm2}) is equivalent to the following proximal step:
\begin{align}
    \theta^{t+1}&= \textbf{prox}_{g/\mt}\left(S^\top w^{t+1}+\tfrac{1}{\mt} \lambda^t\right) \\
    &=\underset{\theta}{\text{argmin}}\left\{g(\theta)+\tfrac{\mt}{2}\norm{S^\top w^{t+1}+\tfrac{1}{\mt}\lambda^t-\theta}^2\right\}.\nonumber
\end{align}
\noindent \textit{Remark}: First note that the Hessian of $\mathcal{L}(\cdot)$ with respect to $w$ can be expressed as: 
\begin{gather}
    H^t=\nabla^2 F(w^t)+\mz D+\mt SS^\top, \label{target hessian}
\end{gather}where $D=A^\top A \in \mathbb{R}^{md}$ is a constant diagonal matrix with its $i$-th block being $\abs{\mathcal{N}_i}I_d$. Moreover, since $F(w)=\sum_{i=1}^m f_i(w_i)$, we conclude that $H^t$ is \textit{block diagonal}. This is only possible by introducing the intermediate consensus variable $z$. If the consensus constraint is directly enforced as $w_i=w_j$, then $H^t$ would have the same sparsity pattern as the network, i.e., the $ij$-th block of $H^t$ would be nonzero if $(i,j)\in \mathcal{E}$. Having a block-diagonal Hessian is \textit{crucial} to eliminate inner-loops for computing quasi-Newton steps. This is because the presence of off-diagonal entries dictates that computing the updates would involve $K$ additional communication rounds among the network, where $K\geq 0$ denotes the number of terms used in the Taylor expansion of $(H^t)^{-1}$ \cite{pdqn}. We proceed to define the $\{s^t_i,q^t_i\}$ pair for L-BFGS approximation used by the $i$-th agent as follows: 
\begin{gather}
q^k_i := \nabla f_i(w^{k+1}_i)-\nabla f_i(w^{k}_i)+\left(\mz\abs{\mathcal{N}_i}+\delta_{il}\mt+\epsilon\right)s_i^k,\nonumber\\  
    s_i^k:=w_i^{k+1}-w_i^k,\label{diff}
\end{gather}
where $\delta_{il}=1$ if $i=l$ and $0$ otherwise, and $\epsilon>0$ brings additional robustness to the algorithm. We emphasize that $q_i^k$ can be computed entirely using local information of agent $i$; this is due to the intermediate consensus variables $\{z_{ij}\}$ as explained above. Moreover, since we approximate the Hessian of the augmented Lagrangian (\ref{target hessian}) instead of just $\nabla^2 F(w^t)$, the update  $(\widehat{H}^t)^{-1}\nabla_w \mathcal{L}$ can be computed at the cost of $\mathcal{O}(cd)$ using Algorithm 1. This is in contradistinction to \cite{pdqn} which opts to approximate $\nabla^2 F(w^t)$ only, whence a linear system has to be solved. We proceed to state a lemma that allows for efficient implementation. 
\begin{algorithm}[t]
	\caption{L-BFGS$-$ADMM} 
	\textbf{Initialization}: zero initialization for all variables.
	\begin{algorithmic}[1]
		\For {$t=0,\dots,T$}
		\hide{\State Agent $i\in [m]$ is active with probability $p_i>0$.}
		\For {each active agent $i$}
		\State Retrieve $w_j,j\in\mathcal{N}_i$ from the buffer.
		\State $h_i \leftarrow \nabla f_i(w_i)+\phi_i+\tfrac{\mz}{2}\sum_{j\in \mathcal{N}_i}\left(w_i-w_j\right)+\delta_{il}\mt(w_i-\theta_i+\tfrac{1}{\mt}\lambda_i)$
		\State \label{step_algo1}Compute the update direction $r_i$ using Algorithm 1 with $h_i$ and  $\{s_i^k,q_i^k\}_{k=\tau_i-c}^{\tau_i-1}$ defined in (\ref{diff}).
		\State \label{primal_step}$w_i\leftarrow w_i-r_i$
		\State \label{phi_step}$\phi_i\leftarrow \phi_i+\tfrac{\mz}{2}\sum_{j\in \mathcal{N}_i}\left(w_i-w_j\right)$
		\State Broadcast $w_i$ to neighbors
		\If {$i=l$} \label{prox_start}
		\State $\theta_i \leftarrow \textbf{prox}_{g_i/\mt}\left(w_i+\tfrac{1}{\mt} \lambda_i\right)$
		\State $\lambda_i \leftarrow \lambda_i+\mt\left(w_i-\theta_i\right)$
		\EndIf \label{prox_end}
		\State Store $\{s_i^{\tau_i},q_i^{\tau_i}\}$ and discard $\{s_i^{\tau_i-c},q_i^{\tau_i-c}\}$
		\EndFor
		\EndFor 
	\end{algorithmic} 
\end{algorithm}

\noindent \textit{Lemma 1}. Define $A_s:=\hat{A}_s\otimes I_d$ $A_d:=\hat{A_d}\otimes I_d$. With zero initialization of $\{y^t,z^t\}$, the following holds: $y^t$ can be decomposed as $y^t=[(\alpha^t)^\top,-(\alpha^t)^\top]^\top$, 
\begin{gather}
    \nabla_w \mathcal{L}^t = \nabla F(w^t)+\phi^t+\tfrac{\mz}{2}L_s w^t+\mt S(S^\top w^t-\theta^t)\nonumber\\+S\lambda^t,\,\,\text{and}\,\,
    z^t = \tfrac{1}{2}(A_s+A_d)w^{t}. \label{lemma1}
\end{gather}
where $\phi^t:=(A_s-A_d)^\top \alpha^t$ and $L_s := (A_s-A_d)^\top (A_s-A_d)$ corresponds to the graph Laplacian. \\
\indent There are two implications from Lemma 1: (i) we only need to store and update half of $y^t$ since it contains entries that are opposite of each other; (ii) there is no need to explicitly store and update $z^{t}$ since it evolves on a manifold defined by $w^t$. 

We explicate the proposed method in Algorithm 2: it admits an asynchronous implementation, where $\tau_i$ represents the counter for agent $i$. Each agent uses a buffer that stores the most recent $\{w_j\}$ values communicated from its neighbors. An active agent $i$ first retrieves $w_j,j\in\mathcal{N}_i,$ and then compute the corresponding subvector $\nabla_w\mathcal{L}^t_i$ in line 4. The quasi-Newton update is then computed in the line \ref{step_algo1} at the cost of $\mathcal{O}(cd)$, where $c$ is the number of $\{s^k_i,q^k_i\}$ pairs stored. After updating $w_i$ (in line \ref{primal_step}), active agent $i$ proceeds to update $\phi_i$ in line \ref{phi_step}. The $l$-th agent additionally performs updates pertaining to the nonsmooth regularization function (lines \ref{prox_start}-\ref{prox_end}). The stored $\{s_i^k,q_i^k\}$ is updated by discarding the one that is computed $\tau_i-c$ steps ago when agent $i$ was activated and adding the most recent $\{s_i^{\tau_i},q_i^{\tau_i}\}$.

\hide{\begin{figure}[htb]
\begin{minipage}[b]{1.0\linewidth}
  \centering
  \centerline{\includegraphics[width=8.5cm]{image1}}
  \centerline{(a) Result 1}\medskip
\end{minipage}
\begin{minipage}[b]{.48\linewidth}
  \centering
  \centerline{\includegraphics[width=4.0cm]{image3}}
  \centerline{(b) Results 3}\medskip
\end{minipage}
\hfill
\begin{minipage}[b]{0.48\linewidth}
  \centering
  \centerline{\includegraphics[width=4.0cm]{image4}}
  \centerline{(c) Result 4}\medskip
\end{minipage}
\caption{Example of placing a figure with experimental results.}
\label{fig:res}
\end{figure}}


\section{Convergence analysis}
\label{sec4}
The analysis is carried assuming a connected network along with the following assumption on local cost functions. 

\noindent \textit{Assumption} 1. Local cost functions $f_i:\mathbb{R}^{d}\to \mathbb{R}, i\in [m],$ are strongly convex with constant $m_f$ and have Lipschitz continuous gradient with constant $M_f$. The local regularizer $g_i(\cdot):\mathbb{R}^d\to \mathbb{R}$ is proper, closed, and convex. 

We denote the unique optimal primal pair as $(w^\star,\theta^\star)$, and the unique dual pair as $(\alpha^\star, \lambda^\star)$, where $\alpha^\star$ lies in the column space of $(A_s-A_d)$. Primal uniqueness follows from strong convexity and we refer the reader to \cite{linearadmm} for existence and uniqueness of dual optimal solutions. We proceed to state a lemma that captures the effect of replacing the exact optimization step (\ref{admm1}) with a quasi-Newton step (\ref{vaniprimal}). 

\noindent \textit{Lemma 2}. Define $L_u:=(A_s+A_d)^\top(A_s+A_d)$ and $E_s:= (A_s-A_d)$. Under Assumption 1, the iterates generated by L-BFGS$-$ADMM and the optimal $(w^\star,\alpha^\star,\lambda^\star)$ satisfy:
\begin{gather}
       \nabla F(w^{t+1})-\nabla F(w^\star)+E_s^\top (\alpha^{t+1}-\alpha^\star)+S\Big(\lambda^{t+1}-\lambda^\star\nonumber\\
       +\mt(\theta^{t+1}-\theta^t)\Big)
       +\left(\tfrac{\mz}{2}L_u+\epsilon I\right)(w^{t+1}-w^t)
       +e^t = 0 \nonumber
\end{gather}
where
$
        e^t = \nabla F(w^t)-\nabla F(w^{t+1})+(\widehat{H}^t-\mz D-\mt SS^\top-\epsilon I) (w^{t+1}-w^t).
$\,\,\textit{Proof}: See the Appendix of \cite{arxiv}. 

Lemma 2 captures the error $e^t$ induced when replacing the primal optimization step (\ref{admm1}) with a quasi-Newton step, i.e., $e^t=0$ corresponds to (\ref{admm1}). We first introduce the following notation before presenting Theorem 1. Denote $u= [w^\top,z^\top,\theta^\top, \alpha^\top, \lambda^\top]^\top$ and similarly $u^\star$ as the optimum, $\sigma_{\mathrm{min}}^{+}$ as the smallest positive eigenvalue of $\begin{bmatrix} E_s \\ S^\top \end{bmatrix}\begin{bmatrix} E_s^\top & S\end{bmatrix}$, and $\sigma_{\mathrm{max}}^G,\sigma_{\mathrm{min}}^G$ as the largest and the smallest eigenvalue of $L_u$. We capture the asynchrony of the algorithm by defining $\Omega^{t}$ as a diagonal random matrix whose blocks activate the corresponding agent (updates $w_i$) or edge (updates $\alpha_i$). We denote $\mathbb{E}^t[\Omega^{t+1}]=\Omega$ and  $p_{\mathrm{min}}:=\underset{i}{\mathrm{min}}\,\,p_i$, with $p_i$ being the probability of each agent and edge being active. \\
\noindent \textit{Theorem 1}. Recall $E_s:= A_s-A_d$ and $L_u:= (A_s+A_d)^\top (A_s+A_d)$ from Lemma 2. Let Assumption 1 hold and assume each agent is activated infinitely often, i,e, $p_{\mathrm{min}}>0$, and let $\mz=2\mt$, the iterates generated by the L-BFGS$-$ADMM using constant initialization $\gamma I$ satisfy:
\begin{gather*}
    \mathbb{E}^t\left[\norm{u^{t+1}-u^\star}_{\mathcal{H}\Omega^{-1}}^2\right] \leq \left(1-\delta\tfrac{p_{\mathrm{min}}}{1+\delta}\right) \norm{u^t-u^\star}_{\mathcal{H}\Omega^{-1}}^2,
\end{gather*}
where $\mathcal{H}:= \textbf{Blkdiag}[\epsilon I,2\mz I,\mt I,\frac{2}{\mz}I,\frac{1}{\mt}I]$, $\tau^t\leq 2\gamma+2c(M_f+\mt(d_{\mathrm{max}}+2)+\epsilon)$, $d_{\mathrm{max}}=\underset{i}{\mathrm{max}}\,\abs{\mathcal{N}_i}$, and $\delta$ satisfy:
  \begin{gather}
       \delta = \min \bigg\{\left(\tfrac{2m_fM_f}{m_f+M_f}-\tfrac{1}{\zeta}\right)\tfrac{1}{\epsilon+\mt(\sigma^{L_u}_{\mathrm{max}}+2)},\tfrac{1}{2},\tfrac{2}{5}\tfrac{\mt\sigma^+_{\mathrm{min}}}{m_f+M_f},\nonumber\\
       \tfrac{\mt\sigma^+_{\mathrm{min}}(\epsilon-\zeta(\tau^t)^2)}{5((\tau^t)^2+\epsilon^2)},\tfrac{\sigma^+_{\mathrm{min}}}{5\max\{1,\sigma^{L_u}_{\max}\}}\bigg\}.
    \end{gather}
\textit{Proof}: See the Appendix of \cite{arxiv}. 
    
\section{Experiments}
\label{sec5}
\begin{figure}[t]
  \centering
  \begin{subfigure}[b]{0.49\linewidth}
    \includegraphics[width=\textwidth]{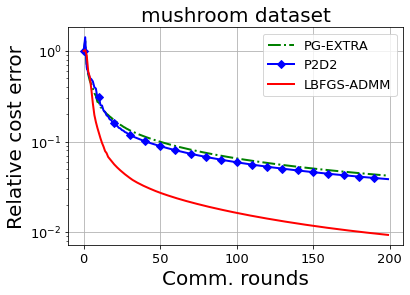}
  \end{subfigure}
  \begin{subfigure}[b]{0.49\linewidth}
    \includegraphics[width=\textwidth]{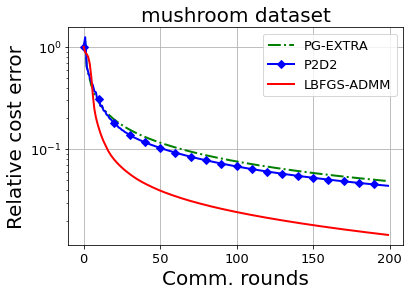}
  \end{subfigure}
  \caption{Performance comparison on the mushrooms dataset with dimension $d=112$. The network size for the left figure is $m=10$ and for the right $m=20$. Each agent of L-BFGS$-$ADMM stores $c=10$ pairs of $\{s_i,q_i\}$. All algorithms are synchronous. }
  \label{fig1}
  \begin{subfigure}[b]{0.49\linewidth}
    \includegraphics[width=\textwidth]{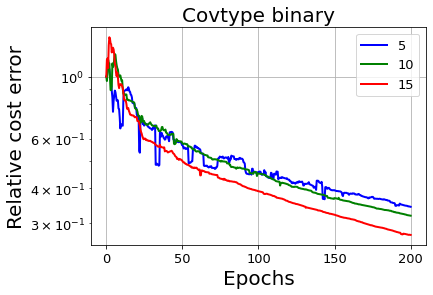}
  \end{subfigure}
  \begin{subfigure}[b]{0.49\linewidth}
    \includegraphics[width=\textwidth]{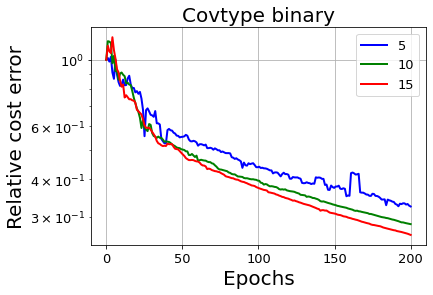}
  \end{subfigure}
  \caption{Effects of number of copies stored and the number of activated agents. We assume a network of 10 agents and activate 2 agents and 5 agents uniformly at random at the left and right plots respectively. We use the same hyperparameters for all cases with $c$ = 5, 10, 15.}
  \label{fig2}
\end{figure}

We evaluate the proposed L-BFGS$-$ADMM against existing distributed algorithms for multi-agent convex composite optimization problems, namely P2D2 \cite{p2d2} and PG-EXTRA \cite{pgextra}. We do not compare against PD-QN \cite{pdqn} since it does not support nonsmooth regularizers and requires $\mathcal{O}(d^3)$ computation costs. We consider the following problem:
\begin{gather*}
    \underset{w\in \mathbb{R}^d}{\text{minimize}}\,\,\mathcal{J}(w)= \left\{\tfrac{1}{m}\sum_{i=1}^m f_i(w)+\rho g(w)\right\},
\end{gather*}
where $
    f_i(w) = \tfrac{1}{m_i}\sum_{j=1}^{m_i} \left[\ln(1+e^{w^\top x_j})+(1-y_j)w^\top x_j\right]$, $g(w) = \norm{w}_1$, and $\rho=0.0005$. We denote $m_i$ as the number of data points held by agent $i$, where each data point contains a feature-label pair $(x_j,y_j)\in \mathbb{R}^d \times\{0,1\}$. We initialize $(\widehat{H}^{t,0})^{-1}$ as a diagonal matrix $\Gamma^t I$, where the $i$-th block is given by 
    $
        \gamma_{i}^t = \frac{(s^{t-1}_i)^\top q_i^{t-1}}{(q_i^{t-1})^\top q_i^{t-1}}.
    $
    Such initialization aims to estimate the norm of the Hessian along the most recent update direction \cite{numerical}. We take 5,000 data points from the covtype dataset with dimension $d=54$ and the mushrooms dataset with dimension $d=112$ respectively. Both datasets are available from the UCI Machine Learning Repository. We plot the averaged relative costs error $\mathrm{RE}(t) :=\tfrac{\tfrac{1}{m}\sum_{i=1}^m\mathcal{J}(w_i^t)-\mathcal{J}(w^\star)}{\tfrac{1}{m}\sum_{i=1}^m\mathcal{J}(w_i^0)-\mathcal{J}(w^\star)}$ in all cases. From Fig \ref{fig1}, we observe that the proposed algorithm outperforms both baseline methods using the number of communication rounds as the metric. On the other hand, Fig. \ref{fig2} shows that storing more copies of $\{s_i,q_i\}$ can effectively reduce oscillations and achieve faster convergence rate.

\vfill\pagebreak

\newpage

\bibliographystyle{IEEE.bst}
\bibliography{main.bib}

\section{Appendix}
In this section, we provide the full proofs for our analytical result. To facilitate the proof, we present some identities that will be useful for our subsequent presentation. The following identities connect our constraint matrix to the graph topology in terms of signed/unsigned incidence matrix $E_s/E_u$ and signed/unsigned graph Laplacian $L_s/L_u$. 
\begin{align*}
    E_s&=A_s-A_d, \,\,E_u=A_s+A_d,\\
    L_s&=E_s^\top E_s,\,\,\,\,\,\,\,\, L_u=E_u^\top E_u,\\
    D&=\tfrac{1}{2}(L_s+L_u)=A_s^\top A_s+A_d^\top A_d.
\end{align*}

\noindent \textit{Assumption} 1. Local cost functions $f_i:\mathbb{R}^{d}\to \mathbb{R}, i\in [m],$ are strongly convex with constant $m_f$ and have Lipschitz continuous gradient with constant $M_f$. The local regularizer $g_i(\cdot):\mathbb{R}^d\to \mathbb{R}$ is proper, closed, and convex. 

\textit{Proof of Lemma 1}: The update for $z^{t+1}$ in (\ref{admm3}) can be obtained by solving the following linear system of equations for $z^{t+1}$:
\begin{align}
    B^\top y^t+\mu_z B^\top(Aw^{t+1}-Bz^{t+1})= 0.\label{z_update}
\end{align}
Recall the dual update for $y^{t+1}$ in (\ref{admm4}), 
\begin{align}
    y^{t+1}=y^t+\mz(Aw^{t+1}-Bz^{t+1}). \tag{5d}
\end{align}
By premultiplying (\ref{admm4}) with $B^\top$ on both sides, we obtain $B^\top y^{t+1}=0$ by (\ref{z_update}). Recalling the definition of the matrix $B$ in (\ref{prob3}) and rewriting $y^{t+1}:=[\alpha^{t+1};\beta^{t+1}]$, we obtain $\beta^{t+1}=-\alpha^{t+1}$ for all $t\geq 0$. With zero initialization of $\{y^t,z^t\}$, $\beta^t=-\alpha^t$ for all $t\geq 0$. Therefore the dual update (\ref{admm4}) can be rewritten as:
\begin{align*}
    \alpha^{t+1} &= \alpha^t+\mz (A_s w^{t+1}-z^{t+1})\\
    -\alpha^{t+1} &= -\alpha^t+\mz (A_d w^{t+1}-z^{t+1}).
\end{align*}
After summing and taking the difference of above, we obtain:
\begin{align}
    z^{t+1} &= \tfrac{1}{2}E_u w^{t+1} \label{new_z},\\
    \alpha^{t+1} &= \alpha^t+\tfrac{\mz}{2}E_sw^{t+1}.
\end{align}
With zero initialization of $\{w^t,z^t\}$, the equation (\ref{new_z}) implies that $z^t=\tfrac{1}{2}(A_s+A_d)w^t=\tfrac{1}{2}E_u w^t$ for all $t\geq 0$. From the definition of the augmented Lagrangian, we obtain: 
\begin{align}
    &\nabla_w \mathcal{L}^t_\mu = \nabla F(w^t)+A^\top y^t+S \lambda^t+\mz A^\top (Aw^t-Bz^t)\nonumber\\
    &+\mt S (S^\top w^t-\theta^t).\label{AL_grad}
\end{align}
Since $y^t=[\alpha^t;-\alpha^t]$ as shown previously and $A=[A_s;A_d]$, we can rewrite $A^\top y^t=E_s^\top \alpha^t=\phi^t$. Moreover, $A^\top Bz^t=\tfrac{1}{2}(A_s+A_d)^\top(A_s+A_d)w^t$. After substituting these into (\ref{AL_grad}), we obtain the desired. \QEDB 

The update for the algorithm can then be summarized as follows: 
\begin{align}
    w^{t+1}&= w^t-(\widehat{H}^t)^{-1}\Big(\nabla F(w^t)+\phi^t+S\lambda^t+\tfrac{\mz}{2}L_s w^t\nonumber\\
    &+\mt S(S^\top w^t-\theta^t)\Big), \label{w_update}\\
    \theta^{t+1} & = \textbf{prox}_{g/\mt}(S^\top w^{t+1}+\tfrac{1}{\mt}\lambda^t) ,\label{theta_update}\\
    \phi^{t+1} &=\phi^t+\tfrac{\mz}{2}L_s w^{t+1}, \label{phi_update}\\
    \lambda^{t+1}& = \lambda^t+\mt (w^{t+1}-\theta^{t+1}).\label{lambda_update}
\end{align}

\textit{Proof of Lemma 2}: From the primal update, the following holds:
\begin{align}
  &\nabla F(w^t)+E_s^\top \alpha^t+S\lambda^t+\tfrac{\mu_z}{2}L_s w^t+\mt S (S^\top w^t-\theta^t)\nonumber\\
    +&\widehat{H}^t(w^{t+1}-w^t)= 0,\label{primal_sub}
\end{align}
where $\phi^t=E_s^\top \alpha^t$. Recall the dual update (\ref{phi_update}). The follwoing holds:
\begin{align}
     E_s^\top \alpha^t+\tfrac{\mu_z}{2}L_sw^t = E_s^\top \alpha^{t+1}-\tfrac{\mu_z}{2}L_s(w^{t+1}-w^t). \label{dual_1}
\end{align}
Similarly, using (\ref{theta_update}), we obtain:
\begin{align}
    &S\lambda^{t}+\mt S^\top(Sw^t-\theta^t) \nonumber\\ =&S\lambda^{t+1}-\mt S\big(S^\top(w^{t+1}-w^t)-(\theta^{t+1}-\theta^t)\big). \label{dual_2}
\end{align}
After substituting (\ref{dual_1}) and (\ref{dual_2}) into (\ref{primal_sub}), we obtain: 
\begin{align}
    &\nabla F(w^t)+E_s^\top \alpha^{t+1}-\tfrac{\mu_z}{2}L_s(w^{t+1}-w^t)+S\big(\lambda^{t+1} \nonumber\\
    &-\mt S^\top(w^{t+1}-w^t)
    +\mt (\theta^{t+1}-\theta^t)\big)\nonumber\\
    &+\widehat{H}^t(w^{t+1}-w^t) = 0 \label{lemma2_1}
\end{align}
After adding and subtracting $(\mu_z D+\mt S S^\top+\epsilon I)(w^{t+1}-w^t)$ from (\ref{lemma2_1}), we obtain
\begin{align*}
    &\nabla F(w^t)+E_s^\top \alpha^{t+1}+\tfrac{\mu_z}{2}L_u(w^{t+1}-w^t)+S\lambda^{t+1}\\
    +&\mu_\theta S(\theta^{t+1}-\theta^t)
    +\epsilon(w^{t+1}-w^t)
    \\
    +&\left(H^t-\mu_z D-\mt S S^\top-\epsilon I\right)(w^{t+1}-w^t)=0.
\end{align*}
Recall the definition of $e^t$:
\begin{align}
           &e^t = \nabla F(w^t)+(\widehat{H}^t-\mz D-\epsilon I-\mt S S^\top)(w^{t+1}-w^t)\nonumber\\
           &-\nabla F(w^{t+1}).\label{error}
\end{align}
After subtracting the following KKT condition and using the definition of $e^t$, 
\begin{gather*}
    \nabla F(w^\star)+E_s^\top \alpha^\star+S\lambda^\star = 0,
\end{gather*}
we obtain the desired. \QEDB

We proceed to establish an upper bound for the error term in the following. 

\noindent \textit{Lemma 3}. The error term (\ref{error}) is upper bounded as follows:
\begin{align*}
     \norm{e^t}\leq  \tau^t\norm{w^{t+1}-w^t},
\end{align*}
where $\tau^t = \norm{\widehat{H}^t-\widehat{H}^{t+1}}$. 

\textit{Proof of Lemma 3}: Recall the secant condition:
\begin{align}
    \widehat{H}^{t+1} s^t = q^t,\label{secant}
\end{align}
where the pair $\{s^t,q^t\}$ is expressed as follows:
\begin{align*}
    s^t &= w^{t+1}-w^t, \\
    q^t &= \nabla F(w^{t+1})-\nabla F(w^t)+(\mz D+\epsilon I+\mt SS^\top ) s^t. 
\end{align*}
Therefore, by substituting the definition of $q^t$ into (\ref{secant}) and rearranging, we obtain: 
\begin{align*}
    &(\widehat{H}^{t+1}-\mz D-\epsilon I-\mt SS^\top) s^t 
    = \nabla F(w^{t+1})-\nabla F(w^t).
\end{align*}
Using the above and the definition of $e^t$ in (\ref{error}), we obtain: 
\begin{align*}
    &\nabla F(w^{t+1}) - \nabla F(w^t)\\
    &= (\widehat{H}^t-\mz D-\epsilon I-\mt S S^\top)s^t-e^t\\
    &= (\widehat{H}^{t+1}-\mz D-\epsilon I-\mt SS^\top) s^t. 
\end{align*}
After rearranging, we obtain:
\begin{gather*}
    e^t = (\widehat{H}^t-\widehat{H}^{t+1})(w^{t+1}-w^t).
\end{gather*}
By applying the Cauchy-Schwartz inequality, we obtain the desired. \QEDB

We proceed to establish the following technical lemma that helps us to bound the curvature estimation obtained from the L-BFGS update. 

\noindent \textit{Lemma 4}. Recall the strong convexity parameter $m_f$ and the Lipschitz smooth parameter $M_f$, where $m_f\leq M_f$. Consider the pair $\{s^t,q^t\}$ used in the construction of the L-BFGS update. The following holds: 
\begin{gather*}
    m_f+\mz+\epsilon\leq \tfrac{\norm{q^t}^2}{(q^t)^\top s^t} \leq  M_f+m\mz +\epsilon+\mt
\end{gather*} 

\noindent \textit{Proof of Lemma 4}: We first define 
\begin{align*}
    r(w^t) = \nabla F(w^t)+(\mz D+\epsilon I+\mt SS^\top) w^t. 
\end{align*}
Then $q^t= r(w^{t+1})-r(w^t)$. Moreover, we denote $r(w^t+\tau (w^{t+1}-w^t))$, $\tau\in [0,1]$, as the $r(\cdot)$ evaluated at some point between $[w^t,w^{t+1}]$. The following holds:
\begin{align*}
    \tfrac{\partial r(\cdot)}{\partial \tau}=[\nabla^2 F(\cdot)+(\mz D+\epsilon I+\mt SS^\top)](w^{t+1}-w^t).
\end{align*}
Therefore, we can rewrite $q^t$ as:
\begin{align*}
    q^t 
    =\int_0^1 [\nabla^2 F(\cdot)+(\mz D+\epsilon I+\mt SS^\top)](w^{t+1}-w^t) d\tau. 
\end{align*}
Denoting $G^t = \int_0^1 [\nabla^2 F(\cdot)+(\mz D+\epsilon I+\mt SS^\top)] d\tau $, we have $G^t s^t =q^t$. Note that $G^t$ can be bounded as follows:
\begin{gather*}
    (m_f+\mz +\epsilon) I \preceq G^t \preceq (M_f+m\mz +\epsilon +\mt) I. 
\end{gather*}
Moreover, we can express $\tfrac{\norm{q^t}^2}{(q^t)^\top s^t}$ as follows:
\begin{align}
    \tfrac{\norm{q^t}^2}{(q^t)^\top s^t}&=\tfrac{(s^t)^\top G^t G^t s^t}{(s^t)^\top G^t s^t}\nonumber\\
    &=\tfrac{(u^t)^\top G^t u^t}{(u^t)^\top u^t},\label{qs_rewrite}
\end{align}
where we denote $u^t = (G^t)^{\tfrac{1}{2}}s^t$. Using the bound for $G^t$ and (\ref{qs_rewrite}), we obtain the desired. \QEDB

We proceed to characterize the curvature estimation obtained by using L-BFGS updates.

\noindent \textit{Lemma 5}. Recall the Lipschitz continuous constant $M_f>0$. Consider the curvature estimation obtained by the L-BFGS updates, $\widehat{H}_i^t\in \mathbb{R}^{d\times d}$. The following holds for all $t\geq 0$:
\begin{gather*}
    \norm{\widehat{H}^t}\leq \gamma+c(M_f+m\mz +\epsilon),
\end{gather*}
where $\gamma$ is the constant used for initialization of L-BFGS and $c$ is the number of copies of $\{s^t_i,q^t_i\}$ pairs, $i\in[m],$ used for constructing $\widehat{H}^t_i$. 

\textit{Proof of Lemma 5}: Recall that the L-BFGS construct the curvature estimation at each agent as follows:
\begin{align}
     \widehat{H}^{t,u+1} = \widehat{H}^{t,u}-\frac{\widehat{H}^{t,u}s^v (s^v)^\top \widehat{H}^{t,u} }{(s^{v})^\top \widehat{H}^{t,u} s^v}+\frac{q^v (q^v)^\top }{(s^v)^\top q^v}, \label{H_construct}
\end{align}
where we have suppressed the subscript $i$. Note that $u\in [0,c]$ denotes the stage of Hessian construction at step $t$, and $v$ denotes the time stamp of stored $\{s^t,q^t\}$ copies, i.e., $v\in [t-c,t-1]$. By taking the trace on both sides of (\ref{H_construct}), we obtain
\begin{align}
      &\mathrm{Tr}\left(\widehat{H}^{t,u+1}\right) = \mathrm{Tr}\left(\widehat{H}^{t,u}\right)-\mathrm{Tr}\left(\frac{\widehat{H}^{t,u}s^v (s^v)^\top \widehat{H}^{t,u} }{(s^{v})^\top \widehat{H}^{t,u} s^v}\right)\nonumber\\
      &+\mathrm{Tr}\left(\frac{q^v (q^v)^\top }{(s^v)^\top q^v}\right). \label{trace}
\end{align}
Note that the second term of (\ref{trace}) can be expressed as:
\begin{gather}
        \mathrm{Tr}\left(\frac{\widehat{H}^{t,u}s^v (s^v)^\top \widehat{H}^{t,u} }{(s^{v})^\top \widehat{H}^{t,u} s^v}\right) = \frac{\norm{\widehat{H}^{t,u}s^v}^2}{(s^{v})^\top \widehat{H}^{t,u} s^v}>0, \label{trace_bd_1}
\end{gather}
where we have used the fact that $\mathrm{Tr}(AB)=\mathrm{Tr}(BA)$. Similarly, the third term can be expressed as:
\begin{align}
    \mathrm{Tr}\left(\frac{q^v (q^v)^\top }{(s^v)^\top q^v}\right)  &= \frac{\norm{q^v}^2}{(s^v)^\top q^u}\nonumber\\
    &\leq M_f+m\mz +\epsilon+\mt,\label{trace_bd_2}
\end{align}
where the inequality follows from Lemma 4. Using (\ref{trace_bd_1}) and (\ref{trace_bd_2}), we obtain an upper bound for (\ref{trace}) as:
\begin{align*}
          \mathrm{Tr}\left(\widehat{H}^{t,u+1}\right) \leq  \mathrm{Tr}\left(\widehat{H}^{t,u}\right)+M_f+d_{\mathrm{max}}\mz+\epsilon+\mt ,
\end{align*}
where $d_{\mathrm{max}}=\underset{i}{\mathrm{max}}$ is the maxim degree. Since the curvature estimation at step $t$ is obtained by applying (\ref{H_construct}) iteratively $c$ times with initialization $\gamma I\succ 0$ and $\widehat{H}^{t,u}\succ 0$, we have $\norm{\widehat{H}^t}=\max_i\{\sigma_i\}<\sum \sigma_i=\mathrm{Tr}(\widehat{H}^t)\leq \gamma+c(M_f+d_{\mathrm{max}}\mz+\epsilon +\mt)$.\QEDB 

The following technical lemma is useful for connecting $\lambda$ and $\theta$. 

\textit{Lemma 6}. Consider the dual variable $\lambda^t\in \mathbb{R}^d$ and the primal variable $\theta^t\in \mathbb{R}^d$. The following holds:
\begin{align*}
    (\lambda^{t+1}-\lambda^t)^\top (\theta^{t+1}-\theta^t)\geq 0 .
\end{align*}

\textit{Proof of Lemma 6}: Recall the update for $\theta^{t+1}$:
\begin{align*}
    \theta^{t+1} &= \textbf{prox}_{g/\mt} (S^\top w^{t+1}+\tfrac{1}{\mt}\lambda^t) \\
    &=\underset{\theta}{\text{argmin}}\bigg\{g(\theta)+\tfrac{\mt}{2}\norm{S^\top w^{t+1}+\tfrac{1}{\mt}\lambda^t-\theta^t}^2\bigg\}.
\end{align*}
The optimality condition gives: 
\begin{align*}
    0 &\in \partial g(\theta^{t+1})- \mt S(S^\top w^{t+1}+\tfrac{1}{\mt}\lambda^t-\theta^{t+1}) \\
    &=\partial g(\theta^{t+1})-\lambda^{t+1}.
\end{align*}
Therefore, we obtain
$
    (\lambda^{t+1}-\lambda^t)^\top (\theta^{t+1}-\theta^t)
    \in (\partial g(\theta^{t+1})-\partial g(\theta^t))^\top  (\theta^{t+1}-\theta^t)\geq 0,
$
where the inequality follows from the convexity of $g(\cdot)$. \QEDB

\noindent \textit{Lemma 7}. Consider the dual variable $[\alpha^t;\lambda^t]\in \mathbb{R}^{(n+1)d}$ and their corresponding optimal pair $[\alpha^\star;\lambda^\star]$ that lies in the column space of $C:=[E_s; S^\top]$. Denote $\sigma^+_{\mathrm{min}}$ as the smallest positive eigenvalue of $CC^\top$ and we select $\mz=2\mt$. The following holds:
\begin{align}
    &\sigma^+_{\mathrm{min}} \big(\norm{\alpha^{t+1}-\alpha^\star}^2+\norm{\lambda^{t+1}-\lambda^\star}^2\big) \leq \nonumber\\
    &4 \big(\norm{\nabla F(w^{t+1})-\nabla F(w^\star)}^2+\norm{w^{t+1}-w^\star}^2_{\overline{L}_u^2}\nonumber\\
    &+\mt^2 \norm{\theta^{t+1}-\theta^t}^2+(\tau^t)^2\norm{w^{t+1}-w^t}^2\big). 
\end{align}

\textit{Proof of Lemma 7}: We first show that $\alpha^t;\lambda^t$ stays in the column space of $C$. Recall the dual update:
\begin{align*}
     \begin{bmatrix}
    \alpha^{t+1}\\
    \lambda^{t+1}
    \end{bmatrix}=
    \begin{bmatrix}
    \alpha^{t}\\
    \lambda^{t}
    \end{bmatrix}+
    \begin{bmatrix}
    \tfrac{\mz}{2}E_s\\
    \mt S^\top
    \end{bmatrix}w^{t+1}-
    \begin{bmatrix}
    0 \\ 
    \mt I_d
    \end{bmatrix}\theta^{t+1}.
\end{align*}
We show that the column space of $[0;I_d]$ is a subspace of the column space of $C$. For fixed $\theta\in \mathbb{R}^d$, we construct $x\in\mathbb{R}^{md}$ as $x:=[\theta;\dots;\theta]$, i.e., $x$ in consensus. Then the following holds:
\begin{align*}
    \begin{bmatrix}
    0 \\ I_d
    \end{bmatrix}\theta = \begin{bmatrix}
    0 \\ \theta
    \end{bmatrix} = \begin{bmatrix}
    E_s \\ S^\top 
    \end{bmatrix} x.
\end{align*}
Since the choice of $\theta$ is arbitrary, we have shown that the column space of $[0;I_d]$ is a subspace of $C$. It is not hard to show that there exists a unique $[\alpha^\star;\lambda^\star]$ in the column space of $C$. Therefore, by choosing $\mz = 2\mt$, we have $[\alpha^{t+1}-\alpha^\star;\lambda^{t+1}-\lambda^\star]$ staying in the column space of $C$. From Lemma 2, the following holds:
\begin{align*}
    & E_s^\top (\alpha^{t+1}-\alpha^\star)+S(\lambda^{t+1}-\lambda^\star) =\\
    &-\big\{\nabla F(w^{t+1})-\nabla F(w^\star)+\overline{L}_u(w^{t+1}-w^t)\nonumber\\
    &+\mt S(\theta^{t+1}-\theta^t)+e^t\big\}
\end{align*}
Since $[\alpha^{t+1}-\alpha^\star;\lambda^{t+1}-\lambda^\star]$ stays in the column space of $C$, $[\alpha^{t+1}-\alpha^\star; \lambda^{t+1}-\lambda^\star]$ is orthogonal to the null space of $C^\top$. Therefore, 
\begin{align*}
     &\sigma^+_{\mathrm{min}} \big(\norm{\alpha^{t+1}-\alpha^\star}^2+\norm{\lambda^{t+1}-\lambda^\star}^2\big) \leq \nonumber\\
    &4 \big(\norm{\nabla F(w^{t+1})-\nabla F(w^\star)}^2+\norm{w^{t+1}-w^\star}^2_{\overline{L}_u^2}\nonumber\\
    &+\mt^2 \norm{\theta^{t+1}-\theta^t}^2+(\tau^t)^2\norm{w^{t+1}-w^t}^2\big),
\end{align*}
where we have used the fact that $(\lambda^{t+1}-\lambda^\star)S^\top S(\lambda^{t+1}-\lambda^\star)=\norm{\lambda^{t+1}-\lambda^\star}^2$, the identity $(\sum_{i=1}^m a_i)^2\leq m \sum_{i=1}^ma_i^2$, and the bound $\norm{e^t}^2\leq (\tau^t)^2\norm{w^{t+1}-w^t}^2$. \QEDB

We proceed to show that the proposed algorithm converges linearly in the following theorem. 

\textit{Proof of Theorem 1}: We first prove the synchronous case. By Assumption 1, the following holds:
\begin{align}
    &\tfrac{m_fM_f}{m_f+M_f}\norm{w^{t+1}-w^\star}^2+\tfrac{1}{m_f+M_f}\norm{\nabla F(w^{t+1})-\nabla F(w^\star)}^2 \nonumber\\
    &\leq (w^{t+1}-w^\star)^\top (\nabla F(w^{t+1})-\nabla F(w^\star)).\label{theorem_1}
\end{align}
Using Lemma 2, we express $\nabla F(w^{t+1})-\nabla F(w^\star)$ as follows:
\begin{align}
    &\nabla F(w^{t+1})-\nabla F(w^\star) =-\big\{E_s^\top (\alpha^{t+1}-\alpha^\star)+S(\lambda^{t+1}-\lambda^\star)\nonumber\\
    &+\overline{L}_u(w^{t+1}-w^t)+\mt S(\theta^{t+1}-\theta^t)+e^t\big\}, \label{theorem_new}
\end{align}
where we have denoted $\overline{L}_u:=\tfrac{\mz}{2}L_u+\epsilon I$. Recall $z^t=\frac{1}{2}E_u w^t$ for all $t\geq 0$ and $L_u=E_u^\top E_u$. Therefore, the following holds:
\begin{align}
    \tfrac{\mz}{2}L_u(w^{t+1}-w^t)= \mz E_u^\top (z^{t+1}-z^t). \label{theorem_new2}
\end{align}By substituting (\ref{theorem_new2}) into (\ref{theorem_new}), and the result into the \textbf{RHS} (right-hand side) of (\ref{theorem_1}), we obtain:
\begin{align}
       &\textbf{RHS}=
        -(w^{t+1}-w^\star)^\top E_s^\top (\alpha^{t+1}-\alpha^\star)-(w^{t+1}-w^\star)^\top e^t&\nonumber\\
       &-(w^{t+1}-w^\star)^\top S(\lambda^{t+1}-\lambda^\star)\nonumber\\
    &-\mt(w^{t+1}-w^\star)^\top S(\theta^{t+1}-\theta^t)\nonumber\\
     &-\epsilon(w^{t+1}-w^\star)^\top (w^{t+1}-w^t)\nonumber\\
     &-\mz(w^{t+1}-w^\star)^\top E_u^\top(z^{t+1}-z^t). \label{theorem_2}
\end{align}
From the dual update, KKT conditions and Lemma 2, the following holds:
\begin{align}
        (w^{t+1}-w^\star)^\top E_s^\top &= \tfrac{2}{\mz}(\alpha^{t+1}-\alpha^t)^\top,\label{theorem_dual_1}\\
    S^\top (w^{t+1}-w^\star) &= \tfrac{1}{\mt}(\lambda^{t+1}-\lambda^t)+\theta^{t+1}-\theta^\star,\label{theorem_dual_2}\\
    \mz (w^{t+1}-w^\star)^\top E_u^\top &=2\mz(z^{t+1}-z^\star)^\top   \label{theorem_dual_3}
\end{align}
After substituting (\ref{theorem_dual_1})-(\ref{theorem_dual_3}) into (\ref{theorem_2}), we obtain:
\begin{align}
        &\textbf{RHS} = -\tfrac{2}{\mz}(\alpha^{t+1}-\alpha^t)^\top (\alpha^{t+1}-\alpha^\star)\nonumber\\
        &-\epsilon (w^{t+1}-w^\star)^\top (w^{t+1}-w^t)\nonumber\\
        &-2\mz(z^{t+1}-z^\star)^\top (z^{t+1}-z^t)\nonumber\\
        &-\tfrac{1}{\mt}(\lambda^{t+1}-\lambda^t)^\top (\lambda^{t+1}-\lambda^\star)\nonumber\\
        &-(\theta^{t+1}
    -\theta^\star)^\top (\lambda^{t+1}-\lambda^\star)-(\lambda^{t+1}-\lambda^t)^\top (\theta^{t+1}-\theta^t)\nonumber\\
    &-\mt(\theta^{t+1}-\theta^\star)^\top (\theta^{t+1}-\theta^t)-(w^{t+1}-w^\star)^\top e^t\nonumber\\
    &\underset{\text{(i)}}{\leq} -\tfrac{2}{\mz}(\alpha^{t+1}-\alpha^t)^\top (\alpha^{t+1}-\alpha^\star)-(w^{t+1}-w^\star)^\top e^t\nonumber\\
    &-\tfrac{1}{\mt}(\lambda^{t+1}-\lambda^t)^\top (\lambda^{t+1}-\lambda^\star)
    \nonumber\\
    &-\mt(\theta^{t+1}-\theta^\star)^\top (\theta^{t+1}-\theta^t) 
    \nonumber\\
    &-\epsilon(w^{t+1}-w^\star)^\top (w^{t+1}-w^t)\nonumber\\
    &-2\mz(z^{t+1}-z^\star)^\top (z^{t+1}-z^t), \label{theorem_3}
\end{align}
where (i) follows from Lemma 6 and the fact that $(\lambda^{t+1}-\lambda^\star)^\top (\theta^{t+1}-\theta^\star)\geq 0$. Using the identity $-2(a-b)^\top  (a-c)=\norm{b-c}^2-\norm{a-b}^2-\norm{a-c}^2$, we can rewrite (\ref{theorem_3}) as:
\begin{align*}
    &2\textbf{RHS} \leq \tfrac{2}{\mz} \big\{\norm{\alpha^t-\alpha^\star}^2-\norm{\alpha^{t+1}-\alpha^\star}^2-\norm{\alpha^{t+1}-\alpha^t}^2\big\}\\
    &+\tfrac{1}{\mt}\big\{\norm{\lambda^t-\lambda^\star}^2-\norm{\lambda^{t+1}-\lambda^\star}^2-\norm{\lambda^{t+1}-\lambda^t}^2\big\}\\
    &+\mt \big\{\norm{\theta^{t}-\theta^\star}^2-\norm{\theta^{t+1}-\theta^\star}^2-\norm{\theta^{t+1}-\theta^t}^2\big\}\\
    &+\epsilon\big\{\norm{w^t-w^\star}^2-\norm{w^{t+1}-w^\star}^2-\norm{w^{t+1}-w^t}^2\big\} \\
    &+2\mz\big\{\norm{z^t-z^\star}^2-\norm{z^{t+1}-z^\star}^2-{z^{t+1}-z^t}^2\big\}\\
    &-2(w^{t+1}-w^\star)^\top e^t \\
    &=\norm{u^t-u^\star}^2_{\mathcal{H}}-\norm{u^{t+1}-u^\star}_{\mathcal{H}}^2-\norm{u^{t+1}-u^t}^2_{\mathcal{H}}\\
    &-2(w^{t+1}-w^t)^\top e^t.
\end{align*}
Recall that \textbf{RHS} stands for the right-hand side of (\ref{theorem_1}). Therefore, the following holds:
\begin{align}
    &\tfrac{2m_fM_f}{m_f+M_f}\norm{w^{t+1}-w^\star}^2+\tfrac{2}{m_f+M_f}\norm{\nabla F(w^{t+1})-\nabla F(w^\star)}^2\nonumber\\
    &+2(w^{t+1}-w^\star)^\top e^t +\norm{u^{t+1}-u^t}^2_{\mathcal{H}}\nonumber\\
    &\leq \norm{u^t-u^\star}^2_{\mathcal{H}}-\norm{u^{t+1}-u^\star}_{\mathcal{H}}^2. \label{theorem_4}
\end{align}
To establish linear convergence, we need to show for some $\delta>0$, the following holds:
\begin{align}
    \delta\norm{u^{t+1}-u^\star}^2_{\mathcal{H}}\leq \norm{u^t-u^\star}^2_{\mathcal{H}}-\norm{u^{t+1}-u^\star}^2_{\mathcal{H}}.\label{linear}
\end{align}
In light of (\ref{theorem_4}), it suffices to show for some $\delta>0$, the following holds:
\begin{align}
    &\delta \norm{u^{t+1}-u^\star}^2_{\mathcal{H}}\leq 
    \tfrac{2m_fM_f}{m_f+M_f}\norm{w^{t+1}-w^\star}^2\nonumber\\
    &+\tfrac{2}{m_f+M_f}\norm{\nabla F(w^{t+1})-\nabla F(w^\star)}^2+\norm{u^{t+1}-u^t}^2_{\mathcal{H}}
    \nonumber\\
    &+2(w^{t+1}-w^\star)^\top e^t.\label{suffice}
\end{align}
We proceed to establish such a bound. Expanding $\norm{u^{t+1}-u^\star}^2_{\mathcal{H}}$, we obtain: 
\begin{align}
    &\norm{u^{t+1}-u^\star}^2_{\mathcal{H}} = \epsilon\norm{w^{t+1}-w^\star}^2+2\mz\norm{z^{t+1}-z^\star}^2\label{u_expand}\\
    &+\mt \norm{\theta^{t+1}-\theta^\star}^2
    +\tfrac{2}{\mz}\norm{\alpha^{t+1}-\alpha^\star}^2+\tfrac{1}{\mt}\norm{\lambda^{t+1}-\lambda^\star}^2\nonumber.  
\end{align}
Note that the last four terms of above are not present in the right-hand of (\ref{suffice}). We proceed to bound these four terms in terms of the components of right-hand side of (\ref{suffice}). Note that by choosing $\mz= 2\mt$, we have $\tfrac{2}{\mz}=\tfrac{1}{\mt}$. Therefore, the following holds: 
\begin{align}
    &\tfrac{2}{\mz}\norm{\alpha^{t+1}-\alpha^\star}^2+\tfrac{1}{\mt}\norm{\lambda^{t+1}-\lambda^\star}^2 \nonumber\\
    &=\tfrac{1}{\mt}\big(\norm{\alpha^{t+1}-\alpha^\star}^2+\norm{\lambda^{t+1}-\lambda^\star}^2\big) \nonumber\\
    &\leq \tfrac{5}{\mt\sigma^+_{\mathrm{min}}}\big(\norm{\nabla F(w^{t+1})-\nabla F(w^\star)}^2+\epsilon^2\norm{w^{t+1}-w^t}^2\nonumber\\
    &+\mt^2 \norm{\theta^{t+1}-\theta^t}^2+(\tau^t)^2\norm{w^{t+1}-w^t}^2\nonumber\\
    &+4\mt^2\sigma^{L_u}_{\mathrm{max}}\norm{z^{t+1}-z^t}^2\big), \label{dual_bd_3}
\end{align}
where the last inequality follows from Lemma 7. Moreover, from the dual update, we obtain:
\begin{align*}
        \mt\norm{\theta^{t+1}-\theta^\star}^2\leq 2\mt \norm{w^{t+1}-w^\star}^2+\tfrac{2}{\mt} \norm{\lambda^{t+1}-\lambda^t}^2.
\end{align*}
Finally, using (\ref{lemma1}), we obtain:
\begin{gather*}
    2\mz\norm{z^{t+1}-z^\star}^2 \leq \tfrac{\mz}{2}\sigma^{L_u}_{\mathrm{max}}\norm{w^{t+1}-w^\star}^2.
\end{gather*}
Substituting these upper bounds for $\tfrac{2}{\mz}\norm{\alpha^{t+1}-\alpha^\star}^2+\tfrac{1}{\mt}\norm{\lambda^{t+1}-\lambda^\star}^2$, $\mt\norm{\theta^{t+1}-\theta^\star}^2$, and $2\mz\norm{z^{t+1}-z^\star}^2$ into (\ref{u_expand}), we obtain:
\begin{align*}
    &\norm{u^{t+1}-u^\star}^2_{\mathcal{H}}\leq (\epsilon+2\mt+\tfrac{\mz\sigma^{L_u}_{\mathrm{max}}}{2})\norm{w^{t+1}-w^\star}^2\\
    &+\tfrac{2}{\mt}\norm{\lambda^{t+1}-\lambda^t}^2
    +\tfrac{5}{\mt\sigma^+_{\mathrm{min}}}\big(\norm{\nabla F(w^{t+1})-\nabla F(w^\star)}^2\\
    &+\epsilon^2\norm{w^{t+1}-w^t}^2
    +\mt^2 \norm{\theta^{t+1}-\theta^t}^2+(\tau^t)^2\norm{w^{t+1}-w^t}^2\nonumber\\
    &+4\mt^2\sigma^{L_u}_{\mathrm{max}}\norm{z^{t+1}-z^t}^2\big).
\end{align*}
Therefore, to establish (\ref{suffice}), it is sufficient to show for some $\delta>0$, the following holds for some $\delta>0$:
\begin{align*}
    &\delta \bigg[(\epsilon+2\mt+\tfrac{\mz\sigma_{\mathrm{max}}^{L_u}}{2})\norm{w^{t+1}-w^\star}^2
    +\tfrac{2}{\mt}\norm{\lambda^{t+1}-\lambda^t}^2\\
    &+\tfrac{5}{\mt\sigma^+_{\mathrm{min}}}\big(\norm{\nabla F(w^{t+1})-\nabla F(w^\star)}^2+\epsilon^2\norm{w^{t+1}-w^t}^2\\
    &+\mt^2 \norm{\theta^{t+1}-\theta^t}^2+(\tau^t)^2\norm{w^{t+1}-w^t}^2\nonumber\\
    &+4\mt^2\sigma^{L_u}_{\mathrm{max}}\norm{z^{t+1}-z^t}^2\big)\bigg]+\zeta\tau^2\norm{w^{t+1}-w^t}^2\\
    &\leq 
    \big(\tfrac{2m_fM_f}{m_f+M_f}-\tfrac{1}{\zeta}\big)\norm{w^{t+1}-w^\star}^2 +\epsilon\norm{w^{t+1}-w^t}^2\\
    &+\tfrac{2}{m_f+M_f}\norm{\nabla F(w^{t+1})-\nabla F(w^\star)}^2
   +2\mz\norm{z^{t+1}-z^t}^2\\
   &+\tfrac{2}{\mz}\norm{\alpha^{t+1}-\alpha^t}^2
      +\tfrac{1}{\mt}\norm{\lambda^{t+1}-\lambda^t}^2
       +\mt \norm{\theta^{t+1}-\theta^t}^2
\end{align*}
where we have used $-\zeta \norm{e^t}^2-\tfrac{1}{\zeta}\norm{w^{t+1}-w^\star}^2\leq 2(w^{t+1}-w^\star)^\top e^t$ holds for any $\zeta>0$, and the bound for $\norm{e^t}$ from Lemma 3, i.e., $\norm{e^t}^2\leq (\tau^t)^2\norm{w^{t+1}-w^t}^2$. The above inequality holds by selecting $\delta$ as:
  \begin{gather}
       \delta = \min \bigg\{\left(\tfrac{2m_fM_f}{m_f+M_f}-\tfrac{1}{\zeta}\right)\tfrac{1}{\epsilon+\mt(\sigma^{L_u}_{\mathrm{max}}+2)},\tfrac{1}{2},\tfrac{2}{5}\tfrac{\mt\sigma^+_{\mathrm{min}}}{m_f+M_f},\nonumber\\
       \tfrac{\mt\sigma^+_{\mathrm{min}}(\epsilon-\zeta(\tau^t)^2)}{5((\tau^t)^2+\epsilon^2)},\tfrac{\sigma^+_{\mathrm{min}}}{5\max\{1,\sigma^{L_u}_{\max}\}}\bigg\}.\label{delta}
    \end{gather}
Note that a uniform upper bound for $\tau^t$ can be obtained by considering Lemma 3 and Lemma 5, i.e., $\tau^t\leq 2\gamma+2c(M_f+d_{\mathrm{max}}\mz+\epsilon)$ for all $t\geq 0$. By selecting $\delta$ as in (\ref{delta}), we establish (\ref{linear}), which proves the linear convergence of the synchronous algorithm. For the asynchronous algorithm, we first express the synchronous algorithm as $u^{t+1}=Tu^{t}$ by defining the operator $T:\mathbb{R}^{(m+2n+2)d}\to\mathbb{R}^{(m+2n+2)d}$. Then the asynchronous algorithm can be expressed as:
\begin{align*}
    u^{t+1} = u^t+\Omega^{t+1}(Tu^t-u^t),
\end{align*}
where
\begin{gather}
    \Omega^{t+1}=\begin{bmatrix}
    X^{t+1}  &      0 & 0            &            0 & 0\\
    0       & Y^{t+1} & 0            &            0 & 0 \\
    0       &       0 & Y^{t+1}      &            0 & 0 \\
    0       &       0 &            0 & X^{t+1}_{ll} & 0\\
    0       &       0 &            0 &            0 & X^{t+1}_{ll}
    \end{bmatrix},
\end{gather}
and $X^{t+1}\in\mathbb{R}^{md\times md},Y^{t+1}\in\mathbb{R}^{nd\times nd}$ are diagonal random matrices with each block $X^{t+1}_{ii}\in\mathbb{R}^{d\times d}, i\in [m]$ and $Y^{t+1}_{kk}\in\mathbb{R}^{d\times d},k\in[n]$, taking values $I_d$ or 0. If we distribute $w_i,i\in[m],$ to agents and $(z_k,\alpha_k),k\in[n]$ to edges, then $w^{t+1}_{i}$ is updated if and only if $X^{t+1}_{ii}=I_d$ and $(z^{t+1}_k,\alpha^{t+1}_k)$ is updated if and only if $Y^{t+1}_{ii}=I_d$. We denote $\mathbb{E}^t[\Omega^{t+1}]=\Omega$. The proof proceeds as follows:
\begin{align}
    &\norm{u^{t+1}-u^\star}^2_{\mathcal{H}\Omega^{-1}}\nonumber\\
    &= \norm{u^{t}+\Omega^{t+1}(Tu^t-u^t)-u^\star}^2_{\mathcal{H}\Omega^{-1}} \nonumber \\
    &= \norm{u^t-u^\star}_{\mathcal{H}\Omega^{-1}}^2 +2(u^t-u^\star)^\top\mathcal{H}\Omega^{-1}\Omega^{t+1}(Tu^t-u^t)\nonumber\\ & + (Tu^t-u^t)^\top\Omega^{t+1}\mathcal{H}\Omega^{-1}\Omega^{t+1}(Tu^t-u^t), \label{theorem_linear_asy}
\end{align}
Since $\Omega^{t+1},\Omega^{-1},\mathcal{H}$ are diagonal matrices, they all commute with each other. Moreover, since the sub-blocks of $\Omega^{t+1}$ are either $I_d$ or $0$, $\Omega^{t+1}\Omega^{t+1}=I_d$. After taking $\mathbb{E}^t[\cdot]$ on both sides of (\ref{theorem_linear_asy}), we obtain:
\begin{align}
     &\mathbb{E}^t\left[\norm{u^{t+1}-u^\star}^2_{\mathcal{H}\Omega^{-1}} \right] \nonumber\\
     &= \norm{u^{t}-u^\star}_{\mathcal{H}\Omega^{-1}}^2
+2(u^t-u^\star)^\top\mathcal{H}(Tu^t-u^t)\nonumber \nonumber\\
&\quad +\norm{Tu^t-u^t}^2_{\mathcal{H}} \nonumber\\ 
&\underset{\text{(i)}}{\leq}  \norm{u^t-u^\star}_{\mathcal{H}\Omega^{-1}}^2-\tfrac{\delta}{1+\delta}\norm{u^t-u^\star}^2_{\mathcal{H}} \nonumber \\
&\underset{\text{(ii)}}{\leq}  \left(1-\tfrac{p^{\mathrm{min}}\delta }{1+\delta} \right)\norm{u^t-u^\star}^2_{\mathcal{H}\Omega^{-1}},
\end{align}
where (i) follows from the fact that $ 2(u-u^\star)^\top\mathcal{H}(Tu-u)+\norm{Tu-u}_{\mathcal{H}}^2 \leq  -\tfrac{\delta}{1+\delta}\norm{u-u^\star}^2_{\mathcal{H}}$ for any $u\in\mathbb{R}^{(m+n+2)d}$ using (\ref{linear}); (ii) $\tfrac{\delta}{1+\delta}\norm{u^t-u^\star}_{\mathcal{H}}\geq \tfrac{p_{\mathrm{min}}\delta}{1+\delta}\norm{u^t-u^\star}^2_{\mathcal{H}\Omega^{-1}}$.   \QEDB

\end{document}